\documentclass[11pt]{article}
\usepackage[latin1]{inputenc}
\usepackage[T1]{fontenc}
\usepackage[english]{babel}
\usepackage{amsmath}
\usepackage{amsthm}
\usepackage{amsfonts}
\usepackage{amssymb}
\usepackage{graphicx}
\usepackage{stmaryrd}
\usepackage{tikz}
\usepackage{url}
\usepackage[np]{numprint}
\usepackage{vmargin}
\usepackage{mathrsfs}
\newtheorem{thm}{Theorem}
\newtheorem{lem}{Lemma}
\newtheorem{prop}{Proposition}
\newtheorem{coro}{Corollary}
\theoremstyle{remark}
\newtheorem{rem}{Remark}
\newcommand{\N}{\mathbb{N}}
\newcommand{\Q}{\mathbb{Q}}
\newcommand{\R}{\mathbb{R}}

\newcommand{\intervalleff }[2]{\left[{#1}\mathpunct{};{#2}\right]}

\newcommand{\intervallefo }[2]{\left[{#1}\mathpunct{};{#2}\right[}
\newcommand{\intervalleoo }[2]{\left]{#1}\mathpunct{};{#2}\right[}

\newcommand{\couple}[2]{\left({#1}\mathpunct{};{#2}\right)}
\newcommand{\triplet}[3]{\left({#1}\mathpunct{};{#2}\mathpunct{};{#3}\right)}
\setmarginsrb{2cm}{2cm}{2cm}{2.5cm}{0cm}{0cm}{0cm}{1cm}
\newcounter{ex}

\begin{document}

\date{}
\title{\textbf{Two trees enumerating the positive rationals }}
\author{Lionel Ponton \\
%EndAName
{\small \texttt{lionel.ponton@gmail.com}}}
\maketitle

\begin{abstract}
We give two trees allowing to represent all positive rational numbers. These trees can be seen as ternary and quinary analogues of the Calkin-Wilf tree. For each of these two trees, we give recurrence formulas allowing to compute the rational number corresponding to the node $n$. These are analogues of the formulas given by Donald Knuth and Moshe Newman for the Calkin-Wilf tree. Finally, we show that the two sequences we have obtained, together with Calkin-Wilf sequence, are the only ones which satisfy a relation analogue to Newman's relation and enumerate the positive rationals.
\end{abstract}

\bigskip

\bigskip

\section{Introduction}

It is well-known, since Cantor's first works on the theory of cardinality, that the rationals are countable. However, it is not so simple to give an explicit enumeration of all of them. Most of the time (see \cite{Bra05}), one proves that $\Q_+$ is countable by constructing a bijection (or an injection) from $\N^{2}$ to $ \N$, which yields an injection from $\Q_+$ to $\N$, and the conclusion follows from Cantor-Bernstein's theorem.

In 2000, N. Calkin and H. S. Wilf \cite{CW00} have described an elegant explicit enumeration of $\Q_{+}^*$. Its first few terms are 
\begin{equation*}
\dfrac{1}{1},~\dfrac{1}{2},~\dfrac{2}{1},~\dfrac{1}{3},~\dfrac{3}{2},~\dfrac{2}{3},~\dfrac{3}{1},~\dfrac{1}{4},~\dfrac{4}{3},~\dfrac{3}{5},~\dfrac{5}{2},~\dfrac{2}{5},~\dfrac{5}{3},~\dfrac{3}{4},~\dfrac{4}{1},~\dfrac{1}{5},~\dfrac{5}{4},~\dfrac{4}{7},...
\end{equation*}

This sequence, known as Calkin-Wilf sequence, is defined by a binary tree the following way~:

\begin{enumerate}
\item[$\bullet $] the top of the tree is $\frac{1}{1}$;

\item[$\bullet $] the vertex labeled $\frac{a}{b}$ has two children~: the left child labeled $\frac{a}{a+b}$ and the right child labeled $\frac{a+b}{b}$.
\end{enumerate}

This leads to the Calkin-Wilf tree, whose first few rows are~:

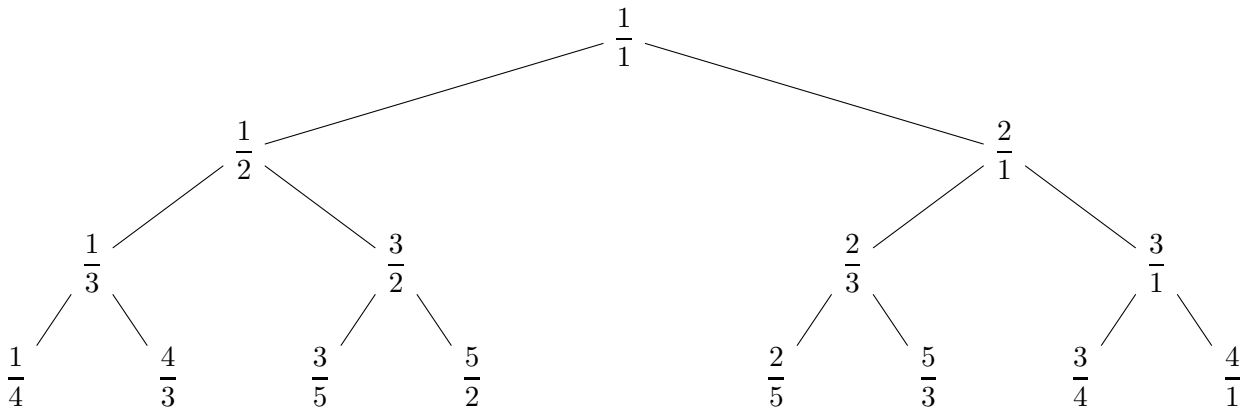
\begin{center}
\begin{tikzpicture}[xscale=1,yscale=1]
\node (P) at ({0},{0}) {$\dfrac{1}{1}$};
\node (fg) at ({-5},{-1.5}) {$\dfrac{1}{2}$};
\node (fgg) at ({-7},{-3}) {$\dfrac{1}{3}$};
\node (fggg) at ({-8},{-4.5}) {$\dfrac{1}{4}$};
\node (fggd) at ({-6},{-4.5}) {$\dfrac{4}{3}$};
\node (fgd) at ({-3},{-3}) {$\dfrac{3}{2}$};
\node (fgdg) at ({-4},{-4.5}) {$\dfrac{3}{5}$};
\node (fgdd) at ({-2},{-4.5}) {$\dfrac{5}{2}$};
\node (fd) at ({5},{-1.5}) {$\dfrac{2}{1}$};
\node (fdg) at ({3},{-3}) {$\dfrac{2}{3}$};
\node (fdgg) at ({2},{-4.5}) {$\dfrac{2}{5}$};
\node (fdgd) at ({4},{-4.5}) {$\dfrac{5}{3}$};
\node (fdd) at ({7},{-3}) {$\dfrac{3}{1}$};
\node (fddg) at ({6},{-4.5}) {$\dfrac{3}{4}$};
\node (fddd) at ({8},{-4.5}) {$\dfrac{4}{1}$};
\draw (P)--(fg);
\draw (fg)--(fgg);
\draw (fgg)--(fggg);
\draw (fgg)--(fggd);
\draw (fg)--(fgd);
\draw (fgd)--(fgdg);
\draw (fgd)--(fgdd);
\draw (P)--(fd);
\draw (fd)--(fdg);
\draw (fdg)--(fdgg);
\draw (fdg)--(fdgd);
\draw (fd)--(fdd);
\draw (fdd)--(fddg);
\draw (fdd)--(fddd);
\end{tikzpicture}
\end{center}

The Calkin-Wilf sequence is then obtained by reading the fraction $\frac{1}{1}$ on level 1, then the two fractions on level 2 from left to right, then the four fractions on level 3 from left to right, and so on. Besides the fact that every positive rational number appears once and only once in reduced form in the tree, this sequence has another remarkable property~: the numerator of the term of rank $n+1$ is equal to the denominator of the term of rank $n$. In other words, there exists a sequence of positive integers $(b_{n})$ such that the term of rank $n$ of the Calkin-Wilf sequence is equal to $\frac{b_{n}}{b_{n+1}}$. In fact, the sequence $(b_{n})$ has been discovered as soon as the mid 19th century, independently by the German mathematician M. Stern \cite{Ste58} and the French clockmaker A. Brocot \cite{Bro61} by considering the median fraction $\frac{a+b}{c+d}$ of two fractions $\frac{a}{b}$ and $\frac{c}{d}$. This procedure leads to another binary tree which enumerates the rationals, named the Stern-Brocot tree \cite[pp. 116-123 et pp. 305-306]{GKP94} and closely connected to the Calkin-Wilf tree (see \cite{Man09} and \cite{BBT10}). B. Reznick \cite{Rez08} notes that Stern has proved in his 1858 paper that, for every pair of positive coprime integers $(a,b)$, there exists one and only one integer $n$ such that $b_{n}=a$ and $b_{n+1}=b$. In other words, Stern proved that $\Q_{+}^*$ is countable more than 15 years before Cantor's first papers on the subject. The sequence $(b_{n})$, which is known nowadays as Stern diatomic sequence, has been widely studied since that time and is known to be connected with many other subjects such as hyperbinary representations, Farey sequences, continued fractions, the
Fibonacci sequence or the Minkowski ?-function (see \cite[pp. 110-114]{AZ06} and \cite{Nor10}).

Calkin-Wilf sequence gives also the answer to a problem set by D. Knuth \cite{Knu01}~: if $v_{p}(n)$ denotes the \textit{p}-adic valuation of the positive integer $n$, prove that the sequence $(x_{n})$ defined by 
\begin{equation}
x_{0}=0\text{ and, for every $n\in \N^*$, }x_{n}=\dfrac{1}{1+2v_{2}(n)-x_{n-1}}  \label{eq22062150}
\end{equation}%
enumerates the positive rationals. Various solutions to this problem have been given in \cite{KRSS03}, among which C. P. Ruppert's one, which associates to the sequence $(x_{n})$ a tree almost identical to Calkin-Wilf tree, the only difference being that the vertices are labeled, not by the
rationals $\frac{a}{b}$, but by the pairs of coprime positive integers $(a,b)$, which is clearly the same. Hence Knuth sequence $(x_{n})$ is exactly the same as Calkin-Wilf sequence.

The editors of \cite{KRSS03} also quote an answer of Moshe Newman, who has shown that the sequence $(x_{n})$ satisfies the recurrence relation~:  
\begin{equation}
\text{For every $n\in \N^*$, }x_{n}=\dfrac{1}{1+2\left\lfloor x_{n-1}\right\rfloor -x_{n-1}}  \label{eq22062207}
\end{equation}%
where $\left\lfloor x\right\rfloor $ denotes the integral part of the real number $x.$ This implies, in particular, the striking result~: 
\begin{equation}
\text{For every $n\in \N^{*}$, }\left\lfloor
x_{n-1}\right\rfloor =v_{2}(n).  
\label{eq22062325}
\end{equation}
Another way to formulate Newman's result consists in saying that the function $f$ defined on $\R_{+}$ by 
\begin{equation}
f:x\mapsto \dfrac{1}{1+2\left\lfloor x\right\rfloor -x}  
\label{eq22062208}
\end{equation}
generates all positive rationals by iteration starting from $x_{0}=0$.

\bigskip

The purpose of this paper is to construct two sequences $(t_{n})$ et $(s_{n}) $ satisfying relations similar to \eqref{eq22062150} and \eqref{eq22062207}. For doing this, we define two trees~: a ternary tree associated to the sequence $(t_{n})$ and a quinary tree associated to the sequence $(s_{n})$. These two trees are not labeled by rationals or pairs of coprime integers, but by triples of integers. They can be considered as generalizations of the Calkin-Wilf tree, in the sense that they lead to sequences which enumerate the postive rationals and satisfy relations similar to \eqref{eq22062150}, \eqref{eq22062207} and \eqref{eq22062325}. However, these generalizations are quite different from those proposed by T. Mansour and M. Shattuck (\cite{MS11} and \cite{MS15}), by B. Bates and T. Mansour \cite{BM11} and by S. H. Chan \cite{Cha11}. Finally, we show that the sequences $(t_{n})$ and $(s_{n}) $ are, together with the Calkin-Wilf
sequence, the only sequences $(u_{n})$ which enumerate the positive rationals and are defined by $u_{0}=0$ and a recurrence relation of the form~: 
\begin{equation}
\text{For every }n\in \N^*,\text{ }u_{n}=\dfrac{f(u_{n-1})}{k},
\label{eq02071820}
\end{equation}
where $f$ is defined by \eqref{eq22062208} and $k\in \N^*$.

\section{A ternary tree}

\subsection{Definition}

We consider the ternary tree $\mathcal{A}_{3}$ whose vertices are labeled by triples of integers $\left( {a}\mathpunct{};{b}\mathpunct{};{c}\right) $ and such that~:

\begin{enumerate}
\item[$\bullet $] the top of the tree is $\left( {1}\mathpunct{};{2}%
\mathpunct{};{0}\right) $;

\item[$\bullet $] the children of $\left( {a}\mathpunct{};{b}\mathpunct{};{c}%
\right) $ are defined by~:

\begin{enumerate}
\item[$\blacktriangleright$] if $b$ is odd~:

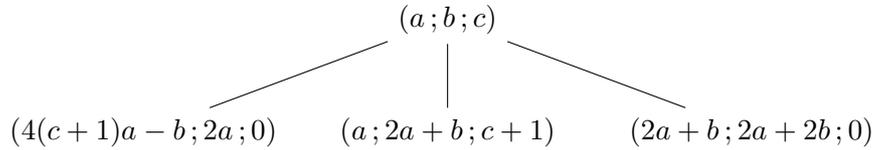
\begin{center}
\begin{tikzpicture}[xscale=1,yscale=1]
\node (P) at ({0},{0}) {$\triplet{a}{b}{c}$};
\node (fg) at ({-4},{-1.5}) {$\triplet{4(c+1)a-b}{2a}{0}$};
\node (fc) at ({0},{-1.5}) {$\triplet{a}{2a+b}{c+1}$};
\node (fd) at ({4},{-1.5}) {$\triplet{2a+b}{2a+2b}{0}$};
\draw (P)--(fg) ;
\draw (P)--(fc) ;
\draw (P)--(fd) ;
\end{tikzpicture}
\end{center}

\item[$\blacktriangleright$] if $b$ is even~:

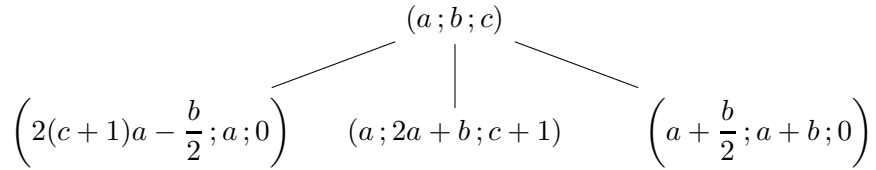
\begin{center}
\begin{tikzpicture}[xscale=1,yscale=1]
\node (P) at ({0},{0}) {$\triplet{a}{b}{c}$};
\node (fg) at ({-4},{-1.5}) {$\triplet{2(c+1)a-\dfrac{b}{2}}{a}{0}$};
\node (fc) at ({0},{-1.5}) {$\triplet{a}{2a+b}{c+1}$};
\node (fd) at ({4},{-1.5}) {$\triplet{a+\dfrac{b}{2}}{a+b}{0}$};
\draw (P)--(fg) ;
\draw (P)--(fc) ;
\draw (P)--(fd) ;
\end{tikzpicture}
\end{center}
\end{enumerate}
\end{enumerate}

\bigskip

Hence the first few levels of $\mathcal{A}_{3}$ are~:

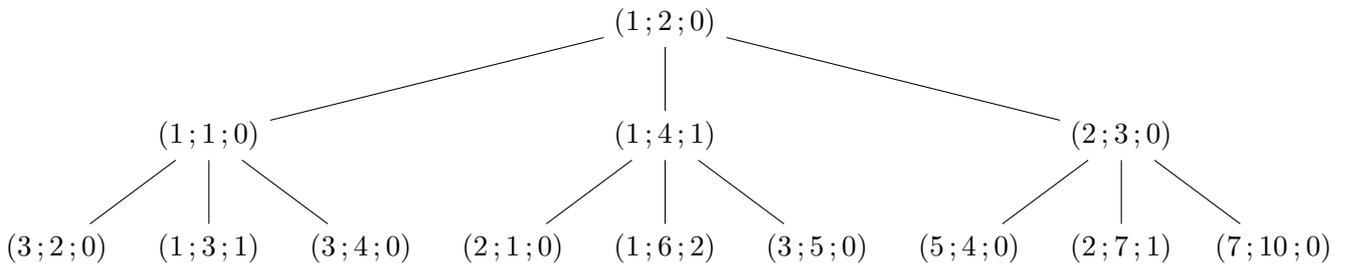
\begin{center}
\begin{tikzpicture}[xscale=1,yscale=1]
\node (P) at ({0},{0}) {$\triplet{1}{2}{0}$};
\node (fg) at ({-6},{-1.5}) {$\triplet{1}{1}{0}$};
\node (fgg) at ({-8},{-3}) {$\triplet{3}{2}{0}$};
\node (fgc) at ({-6},{-3}) {$\triplet{1}{3}{1}$};
\node (fgd) at ({-4},{-3}) {$\triplet{3}{4}{0}$};
\node (fc) at ({0},{-1.5}) {$\triplet{1}{4}{1}$};
\node (fcg) at ({-2},{-3}) {$\triplet{2}{1}{0}$};
\node (fcc) at ({0},{-3}) {$\triplet{1}{6}{2}$};
\node (fcd) at ({2},{-3}) {$\triplet{3}{5}{0}$};
\node (fd) at ({6},{-1.5}) {$\triplet{2}{3}{0}$};
\node (fdg) at ({4},{-3}) {$\triplet{5}{4}{0}$};
\node (fdc) at ({6},{-3}) {$\triplet{2}{7}{1}$};
\node (fdd) at ({8},{-3}) {$\triplet{7}{10}{0}$};
\draw (P)--(fg);
\draw (fg)--(fgg);
\draw (fg)--(fgc);
\draw (fg)--(fgd);
\draw (P)--(fc);
\draw (fc)--(fcg);
\draw (fc)--(fcc);
\draw (fc)--(fcd);
\draw (P)--(fd);
\draw (fd)--(fdg);
\draw (fd)--(fdc);
\draw (fd)--(fdd);
\end{tikzpicture}
\end{center}

The three children of the vertex $N=\left( {a}\mathpunct{};{b}\mathpunct{};{c}\right) $ are called respectively the left, the middle and the right child of $N$ and we say that $N$ is the parent of these three children.

For every $n\in \N^*$, we denote by $N_{n}=\left( {a_{n}}\mathpunct{};{b_{n}}\mathpunct{};{c_{n}}\right) $ the vertex of index $n$ of the tree $\mathcal{A}_{3}$ read from the top and, at each level, from left to right. Hence, $N_{1}=\left( {1}\mathpunct{};{2}\mathpunct{};{0}\right) $, $N_{2}=\left( {1}\mathpunct{};{1}\mathpunct{};{0}\right) $, $N_{3}=\left( {1}\mathpunct{};{4}\mathpunct{};{1}\right) $, and so on... Observe that, by definition, for every $n\in \N^*$, the left, middle and right children of $N_{n}$ are respectively $N_{3n-1}$, $N_{3n}$ and $N_{3n+1}$.

\begin{lem} --- For every $n\in \N^*$, $c_{n}=v_{3}(n)$. 
\label{lem16061603}
\end{lem}

\emph{Proof}. --- For $n=1,$ it is true since $c_{1}=0=v_{3}(1)$. Assume that 
$c_{n}=v_{3}(n)$ for a given $n\in \N^*$. Then the
left child of $N_{n}$ is $N_{3n-1},$ whence by definition $%
c_{3n-1}=0=v_{3}(3n-1)$. Similarly, $N_{3n+1}$ is the right child of $N_{n}$
and $c_{3n+1}=0=v_{3}(3n+1)$. Finally, as $N_{3n}$ is the middle child of $%
N_{n}$, $c_{3n}=c_{n}+1=v_{3}(n)+1=v_{3}(3n)$, and lemma 1 is proved by
induction.\hfill $\square $

\begin{lem}
--- For every $n\in \N^*$, $a_{n}$ and $b_{n}$ are
positive coprime integers and $2a_{n}\geqslant b_{n}-4a_{n}c_{n}$. \label%
{lem16061604}
\end{lem}

\emph{Proof}. --- For $n=1,$ it is true since $a_{1}=1$ , $b_{1}=2$ and $%
c_{1}=0$. Assume that, for a given $n\in \N^*$, $a_{n}$
and $b_{n}$ are positive coprime integers satisfying $2a_{n}\geqslant
b_{n}-4a_{n}c_{n}$.

Assume that $b_{n}$ is odd. Then the three children of $N_{n}$ are
\begin{align*}
N_{3n-1}&=\triplet{4(c_{n}+1)a_{n}-b_{n}}{2a_{n}}{0} \\
N_{3n}&=\triplet{a_{n}}{2a_{n}+b_{n}}{c_{n}+1} \\
N_{3n+1}&=\triplet{2a_{n}+b_{n}}{2a_{n}+2b_{n}}{0}
\end{align*}

As $a_{n}$ and $b_{n}$ are positive integers, it is clear that $b_{3n-1}=2a_{n}$, $a_{3n}=a_{n}$, $b_{3n}=a_{3n+1}=2a_{n}+b_{n}$ and $b_{3n+1}=2a_{n}+2b_{n}$ are positive integers. Moreover, since $2a_{n}\geqslant b_{n}-4a_{n}c_{n}$, 
\begin{equation}
a_{3n-1}=4(c_{n}+1)a_{n}-b_{n}=4a_{n}-(b_{n}-4a_{n}c_{n})\geqslant 2a_{n}
\label{eq16061349}
\end{equation}

and $a_{3n-1}$ is also a positive integer.

Let $d=\text{gcd}\left( a_{3n-1},b_{3n-1}\right) $. Then $d$ divides $b_{3n-1}=2a_{n}$ and $2(c_{n}+1)b_{3n-1}-a_{3n-1}=b_{n}$. Hence $d$ is odd since $b_{n}$ is odd and therefore $d$ divides $a_{n}$. As $a_{n}$ and $b_{n}$ are coprime, we have $d=1$, which means that $\left( a_{3n-1},b_{3n-1}\right) $ are coprime. Similarly we obtain $\left( a_{3n},b_{3n}\right) =\left(a_{3n+1},b_{3n+1}\right) =1.$

Finally, for $N_{3n-1}$ we have, by using (\ref{eq16061349}), 
\begin{equation*}
b_{3n-1}-4a_{3n-1}c_{3n-1}=2a_{n}\leqslant a_{3n-1}\leqslant 2a_{3n-1}.
\end{equation*}

For $N_{3n}$, by using(\ref{eq16061349}), 
\begin{equation*}
b_{3n}-4a_{3n}c_{3n}=b_{n}-4a_{n}c_{n}-2a_{n}\leqslant
2a_{n}-2a_{n}=0\leqslant 2a_{3n},
\end{equation*}

And for $N_{3n+1}$, 
\begin{equation*}
b_{3n+1}-4a_{3n+1}c_{3n+1}=2a_{n}+2b_{n}\leqslant 4a_{n}+2b_{n}=2a_{3n+1}.
\end{equation*}

In the case where $b_{n}$ is even, the proof is similar. One only has to
replace (\ref{eq16061349}) by 
\begin{equation}
a_{3n-1}=2(c_{n}+1)a_{n}-\dfrac{b_{n}}{2}=2a_{n}-\dfrac{b_{n}-4a_{n}c_{n}}{2}%
\geqslant a_{n}.  \label{eq16051429}
\end{equation}

Hence Lemma 2 is proved by induction.\hfill $\square $

\bigskip

Now we put, for every $n\in \N^*$, 
\begin{equation*}
t_{n}=\frac{a_{n}}{b_{n}}.
\end{equation*}

By Lemma 2, $(t_{n})_{n\in \N^*}$ is a sequence of positive reduced rationals. The first few terms of this sequence are~:

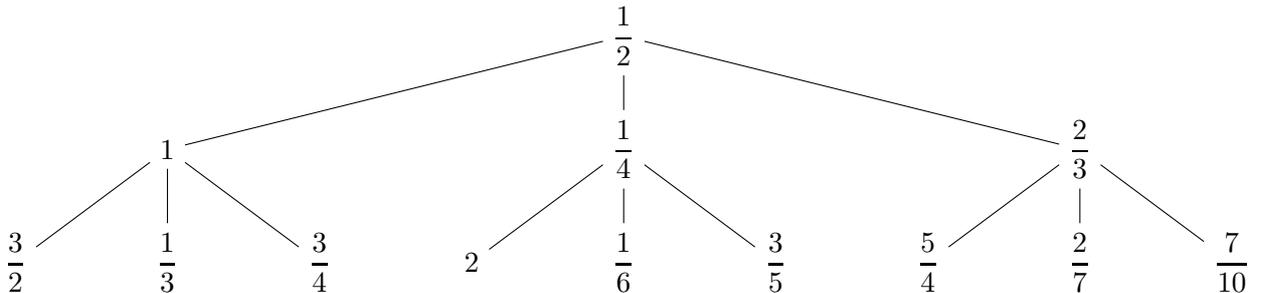
\begin{center}
\begin{tikzpicture}[xscale=1,yscale=1]
\node (P) at ({0},{0}) {$\dfrac{1}{2}$};
\node (fg) at ({-6},{-1.5}) {$1$};
\node (fgg) at ({-8},{-3}) {$\dfrac{3}{2}$};
\node (fgc) at ({-6},{-3}) {$\dfrac{1}{3}$};
\node (fgd) at ({-4},{-3}) {$\dfrac{3}{4}$};
\node (fc) at ({0},{-1.5}) {$\dfrac{1}{4}$};
\node (fcg) at ({-2},{-3}) {$2$};
\node (fcc) at ({0},{-3}) {$\dfrac{1}{6}$};
\node (fcd) at ({2},{-3}) {$\dfrac{3}{5}$};
\node (fd) at ({6},{-1.5}) {$\dfrac{2}{3}$};
\node (fdg) at ({4},{-3}) {$\dfrac{5}{4}$};
\node (fdc) at ({6},{-3}) {$\dfrac{2}{7}$};
\node (fdd) at ({8},{-3}) {$\dfrac{7}{10}$};
\draw (P)--(fg);
\draw (fg)--(fgg);
\draw (fg)--(fgc);
\draw (fg)--(fgd);
\draw (P)--(fc);
\draw (fc)--(fcg);
\draw (fc)--(fcc);
\draw (fc)--(fcd);
\draw (P)--(fd);
\draw (fd)--(fdg);
\draw (fd)--(fdc);
\draw (fd)--(fdd);
\end{tikzpicture}
\end{center}

We remark that, for every $k\in \N^*$,

\begin{gather}
t_{3k-1}=\dfrac{4(c_{k}+1)a_{k}-b_{k}}{2a_{k}}=2(v_{3}(k)+1)-\dfrac{1}{2t_{k}%
}=2v_{3}(3k)-\dfrac{1}{2t_{k}},  \label{eq170610461} \\
t_{3k}=\dfrac{a_{k}}{2a_{k}+b_{k}}=\dfrac{t_{k}}{2t_{k}+1},
\label{eq170610462} \\
t_{3k+1}=\dfrac{2a_{k}+b_{k}}{2a_{k}+2b_{k}}=\dfrac{2t_{k}+1}{2t_{k}+2}.
\label{eq170610463}
\end{gather}

We extend this sequence to $\N$ by putting 
\begin{equation*}
t_{0}=0.
\end{equation*}

We will show that $(t_{n})_{n\in \N}$ enumerates the non negative rationals, i.e. that $n\mapsto t_{n}$ is a bijection from $\N$ to $\Q_{+}$. Before this, we will give two recurrence relations satisfied by the sequence $(t_{n})$.

\subsection{Two recurrence relations}

First we prove that the sequence $(t_{n})$ satisfies a recurrence relation similar to \eqref{eq22062150}.

\begin{prop} --- For every $n\in \N^*$,  $\displaystyle t_{n}=\dfrac{1}{2(1+2v_{3}(n)-t_{n-1})}$.
\label{prop16061604}
\end{prop}

\emph{Proof}. --- This is true for $n=1$ and $n=2$, since
\[\dfrac{1}{2(1+2v_{3}(1)-t_{0})}=\dfrac{1}{2}=t_{1} \qquad \text{and} \qquad \dfrac{1}{2(1+2v_{3}(2)-t_{1})} =1=t_{2}.\]
Now assume that, for a given integer $n\geqslant 3$, the property is true for every positive integer $j\leqslant n-1$. Denote $N_{k}$ $(k\in\N^*$) the parent of $N_{n}$.

\emph{1st case}. --- If $N_{n}$ is the left child of $N_{k},$ then $n=3k-1$ and $N_{n-1}$ is the right child of $N_{k-1}$. As the property is true when $n=k,$ we have by using  \eqref{eq170610461} 
\begin{equation}
t_{n}=2(v_{3}(k)+1)-\dfrac{1}{2t_{k}}=2(v_{3}(k)+1)-\left[
1+2v_{3}(k)-t_{k-1}\right] =1+t_{k-1}.  \label{eq17061056}
\end{equation}

Moreover, since $N_{n-1}$ is the right child of $N_{k-1}$, $t_{n-1}=\frac{2t_{k-1}+1}{2t_{k-1}+2}$ by \eqref{eq170610463}. But $v_{3}(n)=v_{3}(3k-1)=0$, whence 
\begin{equation*}
\dfrac{1}{2(1+2v_{3}(n)-t_{n-1})}=\dfrac{1}{2\left( 1-\frac{2t_{k-1}+1}{2t_{k-1}+2}\right) }=1+t_{k-1}=t_{n}.
\end{equation*}

\emph{2nd case}. --- If $N_{n}$ is the middle child of $N_{k}$, then $n=3k$ and $N_{n-1}$ is the left child of $N_{k}$. By using \eqref{eq170610461}, we have 
\begin{equation*}
t_{n-1}=t_{3k-1}=2v_{3}(3k)-\dfrac{1}{2t_{k}}=2v_{3}(n)-\dfrac{1}{2t_{k}}.
\end{equation*}

Therefore by using \eqref{eq170610462} we obtain 
\begin{equation*}
\dfrac{1}{2(1+2v_{3}(n)-t_{n-1})}=\dfrac{1}{2+\frac{1}{t_{k}}}=\dfrac{t_{k}}{%
2t_{k}+1}=t_{n}.
\end{equation*}

\emph{3rd case}. --- If $N_{n}$ is the right child of $N_{k},$ then $n=3k+1$ and $N_{n-1}$ is the middle child of $N_{k}$. Hence, by \eqref{eq170610462} and \eqref{eq170610463}, 
\begin{equation*}
t_{n-1}=\dfrac{t_{k}}{2t_{k}+1} \qquad \text{and} \qquad t_{n}=\dfrac{2t_{k}+1}{2t_{k}+2}.
\end{equation*}

Since $v_{3}(n)=v_{3}(3k+1)=0$, we have 
\begin{equation*}
\dfrac{1}{2(1+2v_{3}(n)-t_{n-1})}=\dfrac{1}{2\left( 1-\frac{t_{k}}{2t_{k}+1}%
\right) }=\dfrac{2t_{k}+1}{2t_{k}+2}=t_{n}.
\end{equation*}

Proposition \label{prop16061604 copy(1)} is proved by induction.\hfill $\square $

\begin{coro} --- For every $k\in \N^*$, $t_{3k-1}=1+t_{k-1}$. 
\label{coro17061009}
\end{coro}

\emph{Proof}. --- This is exactly the equality (\ref{eq17061056}).

\begin{coro} --- For every $k\in \N^*$, $t_{3k}\in %
\intervalleoo{0}{\tfrac{1}{2}}$, $t_{3k+1}\in \intervalleoo{\tfrac{1}{2}}{1}$
and $t_{3k+2}\in \intervalleoo{1}{+\infty}$. 
\label{coro16061716}
\end{coro}

\emph{Proof}. --- Let $k\in \N^*$. Since $t_{k}>0$, 
\begin{equation*}
t_{3k}=\dfrac{t_{k}}{2t_{k}+1}\in \intervalleoo{0}{\tfrac{1}{2}} \qquad \text{and} \qquad t_{3k+1}=\dfrac{2t_{k}+1}{2t_{k}+2}\in \intervalleoo{\tfrac{1}{2}}{1}. 
\end{equation*}

Moreover, from corollary \ref{coro17061009}, $t_{3k+2}=t_{3(k+1)-1}=1+t_{k}>1 $. \hfill $\square $

\begin{rem} --- As $t_{0}=0$, $t_{1}=\tfrac{1}{2}$ and $t_{2}=1$, we can also see that
for every $k\in \N$, $t_{3k}\in \intervallefo{0}{\tfrac{1}{2}}$, $t_{3k+1}\in \intervallefo{\tfrac{1}{2}}{1}$ et $t_{3k+2}\in \intervallefo{1}{+\infty}$.
\end{rem}

\bigskip

Now we prove that the sequence $(t_{n})$ satisfies relations similar to \eqref{eq22062207} and \eqref{eq22062325}.

\begin{prop} --- For every $n\in \N^*$, $\left\lfloor t_{n-1}\right\rfloor =v_{3}(n)$. \label{prop16061909}
\end{prop}

\emph{Proof}. --- For $n=1$, clearly $\left\lfloor t_{0}\right\rfloor=\left\lfloor 0\right\rfloor =0=v_{3}(1)$.

Let $n\geqslant 2.$ Assume that, for every positive integer $j\leqslant n-1$, $\left\lfloor t_{j-1}\right\rfloor =v_{3}(j)$ and denote by $N_{k}$ the parent of $N_{n}$ ($k\in \N^*$).

If $N_{n}$ is the left child of $N_{k},$ then $n=3k-1$, whence $v_{3}(n)=v_{3}(3k-1)=0$ and by Remark 1 $\left\lfloor t_{n-1}\right\rfloor=\left\lfloor t_{3(k-1)+1}\right\rfloor =0$.

If $N_{n}$ is the right child of $N_{k},$ then $n=3k+1$, $v_{3}(n)=v_{3}(3k+1)=0$ and $\left\lfloor t_{n-1}\right\rfloor =\left\lfloor t_{3k}\right\rfloor =0$.

If $N_{n}$ is the middle child of $N_{k}$ then $n=3k$ and $N_{n-1}$ is the left child of $N_{k}$. Hence, by Corollary \ref{coro17061009}, $t_{n-1}=1+t_{k-1}$ and $\left\lfloor t_{n-1}\right\rfloor =1+\left\lfloor t_{k-1}\right\rfloor $. By the induction hypothesis, it follows that $\left\lfloor t_{n-1}\right\rfloor =1+v_{3}(k)=v_{3}(3k)=v_{3}(n)$.

Proposition 2 is proved by induction.\hfill $\square $

\bigskip

From \ref{prop16061604} and \ref{prop16061909} we get directly

\begin{coro} --- Let $f$ be defined in \eqref{eq22062208}. Then the sequence $(t_{n})_{n\in \N}$ satisfies $t_{0}=0$ and, for every $n\in %
\N^*$, 
\[t_{n}=\dfrac{1}{2(1+2\left\lfloor t_{n-1}\right\rfloor -t_{n-1})}=\dfrac{f(t_{n-1})}{2}.\]
\end{coro}

\subsection{The sequence $(t_{n})$ enumerates $\Q_{+}$}

\begin{thm} --- The mapping $n\mapsto t_{n}$ is a bijection from $\N$ to $\Q_+$. \label{th16061912}
\end{thm}

\emph{Proof}. --- As $t_{0}=0$ and $t_{n}=\frac{a_{n}}{b_{n}}$ is reduced for every $n\in \N^*$, we have to prove that, for every pair of non zero coprime natural integers $\left( {\alpha }\mathpunct{};{\beta }\right) $, there exists one and only one $n\geqslant 1$ such that $a_{n}=\alpha $ and $b_{n}=\beta $.

The proof is by induction on $m=\alpha +\beta $.

If $m=2$ then $\alpha =\beta =1$ and Corollary \ref{coro16061716} implies that $n=2$ is the only one integer such that $a_{n}=b_{n}=1$.

Assume that, for a given integer $m\geqslant 2$, the property is true for every $k\in \{2,...,m\}$. Let $\left( {\alpha }\mathpunct{};{\beta }\right)$ be a pair of non zero coprime natural integers such that $\alpha +\beta =m+1$.

\emph{1st case}~: $\beta >2\alpha $. Then, by Corollary \ref{coro16061716}, if $n$ exists, there exists $k\in \N^*$ such that $n=3k$. Hence $N_{n}$ is the middle child of $N_{k}$. Therefore $%
N_{k}=\left( {\alpha }\mathpunct{};{\beta -2\alpha }\mathpunct{};{c_{k}}\right) $. Now, $\alpha +(\beta -2\alpha )=\beta -\alpha \leqslant m$ and $\alpha $ and $\beta -2\alpha $ are coprime. By the induction hypothesis, there exists one and only one integer $k$ such that $a_{k}=\alpha $ and $b_{k}=\beta -2\alpha ,$ which proves that $n=3k$ is the one and only one integer such that $a_{n}=\alpha $ and $b_{n}=\beta $.

\emph{2nd case}~: $\beta =2\alpha .$ Then, $\couple{\alpha}{\beta}=\couple{1}{2}$ since $\alpha$ and $\beta$ are coprime. By Corollary \ref{coro16061716}, $n=1$ is the sole integer such that $a_{n}=1$ and $b_{n}=2$. 

\emph{3rd case}~: $\beta <2\alpha <2\beta .$ Then, by Corollary \ref{coro16061716}, if $n$ exists, there exists $k\in \N^*$ such that $n=3k+1$. Hence $N_{n}$ is the right child of $N_{k}$. If $\beta $ is even $N_{k}=\left( {\alpha -\frac{\beta }{2}}\mathpunct{};{\beta -\alpha }\mathpunct{};{c_{k}}\right) $. Since $\alpha -\frac{\beta }{2}+\beta -\alpha=\frac{\beta }{2}\leqslant m$ and $\alpha -\frac{\beta }{2}$ and $\beta-\alpha $ are coprime, we see, as in the first case, that $n=3k+1$ is the one and only one integer such that $a_{n}=\alpha $ and $b_{n}=\beta $. If $\beta $ is odd then $N_{k}=\left( {2\alpha -\beta }\mathpunct{};{2\beta-2\alpha }\mathpunct{};{c_{k}}\right) .$ As $2\alpha -\beta +2\beta -2\alpha=\beta \leqslant m$ and $2\alpha -\beta $ and $2\beta -2\alpha $ are coprime (since $\beta $ is even), we draw the same conclusion.

\emph{4th case}~: $\alpha =\beta .$ Then $\alpha =\beta =1$ since $\alpha $ and $\beta $ are coprime. But this is impossible because $\alpha+\beta =m+1\geqslant 3$.

\emph{5th case}~: $\alpha >\beta .$ Then, by Corollary \ref{coro16061716}, if $n$ exists, there exists an integer $k\geqslant 2$ such that $n=3k-1$. Hence $N_{n}$ is the left child of $N_{k}$. In this case, we cannot argue as before because, for odd $b_{n}$, $a_{n}+b_{n}$ is not necessarily greater
than $a_{k}+b_{k},$ as can be seen, for example, when $N_{3}=\left( {1}\mathpunct{};{4}\mathpunct{};{1}\right) $ and $N_{8}=\left( {2}\mathpunct{};{1}\mathpunct{};{0}\right) $. However, by Corollary \ref{coro17061009}, $t_{n}=1+t_{k-1}$, whence $t_{k-1}=\frac{\alpha -\beta }{\beta }$. As $(\alpha -\beta )+\beta =\alpha \leqslant m$ and $\alpha -\beta $ and $\alpha 
$ are coprime, by the induction hypothesis there exists one and only one integer $k\geqslant 2$ such that $a_{k-1}=\alpha -\beta $ and $b_{k-1}=\beta.$ This shows that $n=3k$ is the one and only one integer such that $a_{n}=\alpha $ and $b_{n}=\beta $.

Theorem 1 is therefore proved by induction.\hfill $\square $

\bigskip

Hence the ternary tree $\mathcal{A}_{3}$ enabled us to construct a sequence $(t_{n})$ which enumerates the non negative rationals and satisfies recurrence relations similar to \eqref{eq22062150} and \eqref{eq22062207}. Now we give a similar construction by using a quinary tree.

\section{A quinary tree}

\subsection{Definition}

We consider the quinary tree $\mathcal{A}_{5}$ whose vertices are labeled by
triples of integers $\left( {a}\mathpunct{};{b}\mathpunct{};{c}\right) $
such that~:

\begin{enumerate}
\item[$\bullet $] the top of the tree is $\left( {1}\mathpunct{};{3}\mathpunct{};{0}\right) $;
\item[$\bullet $] the children of $\left( {a}\mathpunct{};{b}\mathpunct{};{c}\right) $ are defined by~:
\newpage
\begin{itemize}
\item[$\blacktriangleright$] if $3$ doesn't divide $b$~:
\end{itemize}
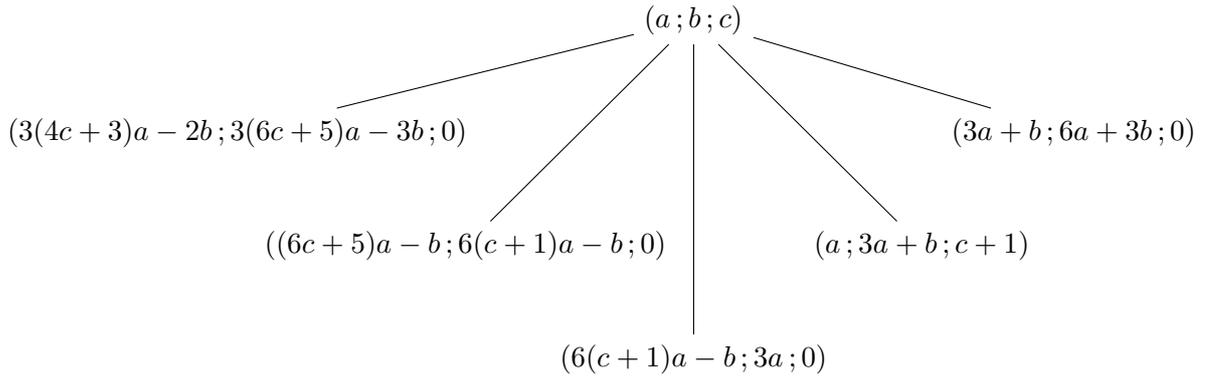
\begin{center}
\begin{tikzpicture}[xscale=1,yscale=1]
\node (P) at ({0},{0}) {$\triplet{a}{b}{c}$};
\node (f1) at ({-6},{-1.5}) {$\triplet{3(4c+3)a-2b}{3(6c+5)a-3b}{0}$};
\node (f2) at ({-3},{-3}) {$\triplet{(6c+5)a-b}{6(c+1)a-b}{0}$};
\node (f3) at ({0},{-4.5}) {$\triplet{6(c+1)a-b}{3a}{0}$};
\node (f4) at ({3},{-3}) {$\triplet{a}{3a+b}{c+1}$};
\node (f5) at ({5},{-1.5}) {$\triplet{3a+b}{6a+3b}{0}$};
\draw (P)--(f1) ;
\draw (P)--(f2) ;
\draw (P)--(f3) ;
\draw (P)--(f4) ;
\draw (P)--(f5) ;
\end{tikzpicture}
\end{center}

\begin{itemize}
\item[$\blacktriangleright$] if 3 divides $b$~:
\end{itemize}

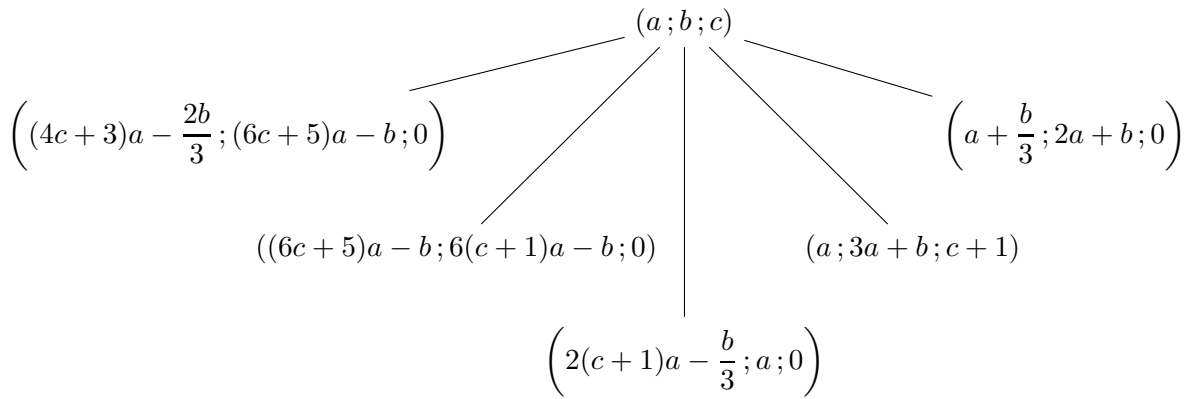
\begin{center}
\begin{tikzpicture}[xscale=1,yscale=1]
\node (P) at ({0},{0}) {$\triplet{a}{b}{c}$};
\node (f1) at ({-6},{-1.5}) {$\triplet{(4c+3)a-\dfrac{2b}{3}}{(6c+5)a-b}{0}$};
\node (f2) at ({-3},{-3}) {$\triplet{(6c+5)a-b}{6(c+1)a-b}{0}$};
\node (f3) at ({0},{-4.5}) {$\triplet{2(c+1)a-\dfrac{b}{3}}{a}{0}$};
\node (f4) at ({3},{-3}) {$\triplet{a}{3a+b}{c+1}$};
\node (f5) at ({5},{-1.5}) {$\triplet{a+\dfrac{b}{3}}{2a+b}{0}$};
\draw (P)--(f1) ;
\draw (P)--(f2) ;
\draw (P)--(f3) ;
\draw (P)--(f4) ;
\draw (P)--(f5) ;
\end{tikzpicture}
\end{center}
\end{enumerate}

The five children of the vertex $N=\left( {a}\mathpunct{};{b}\mathpunct{};{c}\right) $ are called from left to right respectively first, second, third,
fourth and fifth child of $N.$

For every $n\in \N^*$, we denote $N_{n}=\left( {a_{n}}\mathpunct{};{b_{n}}\mathpunct{};{c_{n}}\right) $ the vertex of index $n$ of the tree $\mathcal{A}_{5}$ read from the top and, at each level, from left to right. Thus, $N_{1}=\left( {1}\mathpunct{};{3}\mathpunct{};{0}\right) $, $N_{2}=\left( {1}\mathpunct{};{2}\mathpunct{};{0}\right) $, $N_{3}=\left( {2}\mathpunct{};{3}\mathpunct{};{0}\right) $, $N_{4}=\left( {1}\mathpunct{};{1}\mathpunct{};{0}\right) $, $N_{5}=\left( {1}\mathpunct{};{6}\mathpunct{};{1}\right) $, and so on...

By definition, for every $n\in \N^*$, the \textit{i}-th child of $N_{n}$ is $N_{5(n-1)+i+1}$.

\bigskip

It is easy to check, as in Lemmas \ref{lem16061603} and \ref{lem16061604}, that for every $n\in \N^*$, $c_{n}=v_{5}(n)$, $a_{n}\in \N^*$, $b_{n}\in \N^*$ (with, this time, $3a_{n}\geqslant b_{n}-6a_{n}c_{n})$ and $\text{gcd}(a_{n},b_{n})=1$. Hence, by putting for every $n\in \N^*$, $
s_{n}=\frac{a_{n}}{b_{n}}$, we define a sequence $(s_{n})_{n\in \N^*}$ of positive reduced rationals, whose first few terms are~:

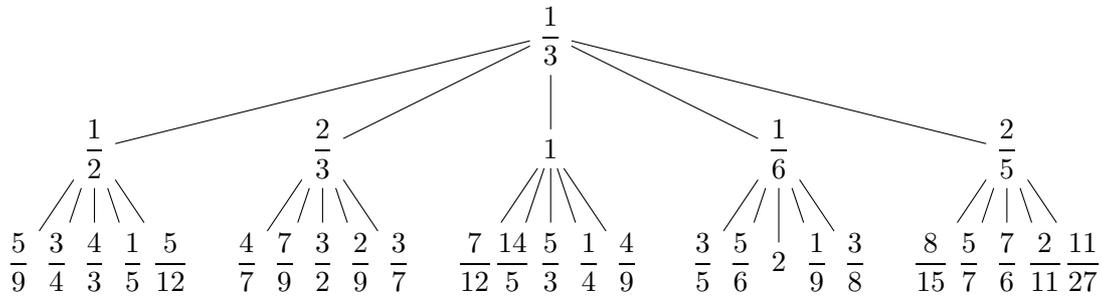
\begin{center}
\begin{tikzpicture}[xscale=1,yscale=1]
\node (P) at ({0},{0}) {$\dfrac{1}{3}$};
\node (f1) at ({-6},{-1.5}) {$\dfrac{1}{2}$};
\node (f11) at ({-7},{-3}) {$\dfrac{5}{9}$};
\node (f12) at ({-6.5},{-3}) {$\dfrac{3}{4}$};
\node (f13) at ({-6},{-3}) {$\dfrac{4}{3}$};
\node (f14) at ({-5.5},{-3}) {$\dfrac{1}{5}$};
\node (f15) at ({-5}, {-3}) {$\dfrac{5}{12}$};
\node (f2) at ({-3},{-1.5}) {$\dfrac{2}{3}$};
\node (f21) at ({-4},{-3}) {$\dfrac{4}{7}$};
\node (f22) at ({-3.5},{-3}) {$\dfrac{7}{9}$};
\node (f23) at ({-3},{-3}) {$\dfrac{3}{2}$};
\node (f24) at ({-2.5},{-3}) {$\dfrac{2}{9}$};
\node (f25) at ({-2}, {-3}) {$\dfrac{3}{7}$};
\node (f3) at ({0},{-1.5}) {$1$};
\node (f31) at ({-1},{-3}) {$\dfrac{7}{12}$};
\node (f32) at ({-0.5},{-3}) {$\dfrac{14}{5}$};
\node (f33) at ({0},{-3}) {$\dfrac{5}{3}$};
\node (f34) at ({0.5},{-3}) {$\dfrac{1}{4}$};
\node (f35) at ({1}, {-3}) {$\dfrac{4}{9}$};
\node (f4) at ({3},{-1.5}) {$\dfrac{1}{6}$};
\node (f41) at ({2},{-3}) {$\dfrac{3}{5}$};
\node (f42) at ({2.5},{-3}) {$\dfrac{5}{6}$};
\node (f43) at ({3},{-3}) {$2$};
\node (f44) at ({3.5},{-3}) {$\dfrac{1}{9}$};
\node (f45) at ({4}, {-3}) {$\dfrac{3}{8}$};
\node (f5) at ({6},{-1.5}) {$\dfrac{2}{5}$};
\node (f51) at ({5},{-3}) {$\dfrac{8}{15}$};
\node (f52) at ({5.5},{-3}) {$\dfrac{5}{7}$};
\node (f53) at ({6},{-3}) {$\dfrac{7}{6}$};
\node (f54) at ({6.5},{-3}) {$\dfrac{2}{11}$};
\node (f55) at ({7}, {-3}) {$\dfrac{11}{27}$};
\draw (P)--(f1);
\draw (f1)--(f11);
\draw (f1)--(f12);
\draw (f1)--(f13);
\draw (f1)--(f14);
\draw (f1)--(f15);
\draw (P)--(f2);
\draw (f2)--(f21);
\draw (f2)--(f22);
\draw (f2)--(f23);
\draw (f2)--(f24);
\draw (f2)--(f25);
\draw (P)--(f3);
\draw (f3)--(f31);
\draw (f3)--(f32);
\draw (f3)--(f33);
\draw (f3)--(f34);
\draw (f3)--(f35);
\draw (P)--(f4);
\draw (f4)--(f41);
\draw (f4)--(f42);
\draw (f4)--(f43);
\draw (f4)--(f44);
\draw (f4)--(f45);
\draw (P)--(f5);
\draw (f5)--(f51);
\draw (f5)--(f52);
\draw (f5)--(f53);
\draw (f5)--(f54);
\draw (f5)--(f55);
\end{tikzpicture}
\end{center}

It can be remarked that, for every $k\in \N^*$,
whether 3 divides $b_{k}$ or not,

\begin{gather}
s_{5k-3}=\dfrac{3(4c_k+3)a_k-2b_k}{3(6c_k+5)a_k-3b_k}=\dfrac{3(4c_k+3)-\frac{%
2}{s_k}}{3(6c_k+5)-\frac{3}{s_k}},  \label{eq200617231} \\
s_{5k-2}=\dfrac{(6c_k+5)a_k-b_k}{6(c_k+1)a_k-b_k}=\dfrac{6c_k+5-\frac{1}{s_k}%
}{6(c_k+1)-\frac{1}{s_k}},  \label{eq200617232} \\
s_{5k-1}=\dfrac{6(c_k+1)a_k-b_k}{3a_k}=2(c_k+1)-\dfrac{1}{3s_k},
\label{eq200617233} \\
s_{5k}=\dfrac{a_k}{3a_k+b_k}=\dfrac{1}{3+\frac{1}{s_k}},  \label{eq200617234}
\\
\quad s_{5k+1}=\dfrac{3a_k+b_k}{6a_k+3b_k}=\dfrac{3s_k+1}{6s_k+3}=\dfrac{3+%
\frac{1}{s_k}}{6+\frac{3}{s_k}}.  \label{eq200617235}
\end{gather}

\bigskip

We extend this sequence to $\N$ by putting $s_{0}=0$. We
will now show, as we did for $(t_{n})_{n\in\N}$ that $(s_{n})_{n\in\N}$, enumerates the elements of $\Q_{+}$.

\subsection{Recurrence relations}

\begin{prop} --- For every $n\in \N^*$, $\displaystyle s_{n}=\dfrac{1}{3(1+2v_{5}(n)-s_{n-1})}.$
\label{prop200954}
\end{prop}

\emph{Proof}. --- For $n=1$ and $n=2,$ this is true since 
\[\dfrac{1}{3(1+2v_{5}(1)-s_{0})}=\dfrac{1}{3}=s_{1} \qquad \text{et} \qquad \dfrac{1}{3(1+v_{5}(2)-s_{1})} =\dfrac{1}{2}=s_{2}.\]
Assume that, for a given $n\geqslant 2$, the property is true for every positive integer $j\leqslant n-1$. Denote $N_{k}$ ($k\in \N^*$) the
parent of $N_{n}$ .

\emph{1st case}. --- If $N_{n}$ is the first child of $N_{k}$ then $n=5k-3$ and $N_{n-1}$ is the fifth child of $N_{k-1}$. By the induction hypothesis, 
\begin{equation*}
s_{k}=\dfrac{1}{3(1+2v_{5}(k)-s_{k-1})}=\dfrac{1}{3\left(
1+2c_{k}-s_{k-1}\right) }.
\end{equation*}%
Hence, by using \eqref{eq200617231}, 
\begin{equation}
s_{n}=\dfrac{3(4c_{k}+3)-\frac{2}{s_{k}}}{3(6c_{k}+5)-\frac{3}{s_{k}}}=\dfrac{3(4c_{k}+3)-6(2c_{k}+1-s_{k-1})}{3(6c_{k}+5)-9(2c_{k}+1-s_{k-1})}=\dfrac{1+2s_{k-1}}{2+3s_{k-1}}.  
\label{eq21060933}
\end{equation}
As $v_{5}(n)=0$, \eqref{eq200617235} yields 
\begin{equation*}
\dfrac{1}{3(2v_{5}(n)+1-s_{n-1})}=\dfrac{1}{3(1-s_{5(k-1)+1})}=\dfrac{1}{3\left( 1-\frac{3s_{k-1}+1}{6s_{k-1}+3}\right) }=\dfrac{2s_{k-1}+3}{3s_{k-1}+2}=s_{n}.
\end{equation*}

\emph{2nd case}. --- If $N_{n}$ is the second child of $N_{k}$ then $%
n=5k-2$ and $N_{n-1}$ is the first child of $N_{k}$. As $v_{5}(n)=0$, %
\eqref{eq200617231} and \eqref{eq200617232} yield 
\begin{equation*}
\dfrac{1}{3(2v_{5}(n)+1-s_{n-1})}=\dfrac{1}{3(1-s_{5k-3})}=\dfrac{1}{3\left(
1-\frac{3(4c_{k}+3)a_{k}-2b_{k}}{3(6k+5)a_{k}-3b_{k}}\right) }=\dfrac{%
(6c_{k}+5)a_{k}-b_{k}}{6(c_{k}+1)a_{k}-b_{k}}=s_{5k-2}=s_{n}.
\end{equation*}

\emph{3rd case}. --- If $N_{n}$ is the third child of $N_{k}$ then $%
n=5k-1$ and $N_{n-1}$ is the second child of $N_{k}$. As $v_{5}(n)=0$, %
\eqref{eq200617232} and \eqref{eq200617233} yield 
\begin{equation*}
\dfrac{1}{3(2v_{5}(n)+1-s_{n-1})}=\dfrac{1}{3(1-s_{5k-2})}=\dfrac{1}{3\left(
1-\frac{(6c_{k}+5)a_{k}-b_{k}}{6(c_{k}+1)a_{k}-b_{k}}\right) }=\dfrac{%
6(c_{k}+1)a_{k}-b_{k}}{3a_{k}}=s_{5k-1}=s_{n}.
\end{equation*}

\emph{4th case}. --- If $N_{n}$ is the fourth child of $N_{k}$ then $%
n=5k$ and $N_{n-1}$ is the third child of $N_{k}$. As $%
v_{5}(5k)=v_{5}(k)+1=c_{k}+1$, \eqref{eq200617233} and \eqref{eq200617234}
yield 
\begin{equation*}
\dfrac{1}{3(2v_{5}(n)+1-s_{n-1})}=\dfrac{1}{3\left( 2(c_{k}+1)+1-\left[
2(c_{k}+1)-\frac{1}{3s_{k}}\right] \right) }=\dfrac{1}{3+\frac{1}{s_{k}}}%
=s_{n}.
\end{equation*}

\emph{5th cas}. --- If $N_{n}$ is the fifth child of $N_{k}$ then $%
n=5k+1$ and $N_{n-1}$ is the fourth child of $N_{k}$. As $v_{5}(n)=0$, %
\eqref{eq200617234} and \eqref{eq200617235} yield 
\begin{equation*}
\dfrac{1}{3(2v_{5}(n)+1-s_{n-1})}=\dfrac{1}{3\left( 1-\frac{a_{k}}{%
3a_{k}+b_{k}}\right) }=\dfrac{3a_{k}+b_{k}}{6a_{k}+3b_{k}}=s_{n}.
\end{equation*}

Proposition 3 is proved by induction.\hfill $\square $

\begin{coro}
--- For every $k\in \N^*$, $s_{5k-1}=1+s_{k-1}$. \label%
{coro20061836}
\end{coro}

\emph{Proof}. --- Let $k$ be a positive integer. By definition, $s_{5k-1}=2(c_{k}+1)-\frac{1}{3s_{k}}$ and, by Proposition \ref{prop200954}, $s_{k}=\frac{1}{3(2c_{k}+1-s_{k-1})}$. Therefore $s_{5k-1}=2(c_{k}+1)-(2c_{k}+1-s_{k-1})=1+s_{k-1}$. \hfill $\square $

\bigskip

As in Corollary \ref{coro16061716}, we deduce from Proposition \ref{prop200954}
that the rationals $s_{n}$ belong to one of the five intervals $\intervallefo{\tfrac{1}{2}}{\tfrac{2}{3}}$, $\intervallefo{\tfrac{2}{3}}{1}$, $\intervallefo{1}{+\infty}$, $\intervallefo{0}{\tfrac{1}{3}}$ or $\intervallefo{\tfrac{1}{3}}{\tfrac{1}{2}}$ depending on their rank in the tree $\mathcal{A}_{5}$ as a first, second, third, fourth or fifth child~:

\begin{coro} --- For every $k\in \N^*$, $s_{5k}\in \intervalleoo{0}{\tfrac{1}{3}}$, $s_{5k+1}\in \intervalleoo{\tfrac{1}{3}}{\tfrac{1}{2}}$, $s_{5k+2}\in \intervalleoo{\tfrac{1}{2}}{\tfrac{2}{3}}$, $s_{5k+3}\in \intervalleoo{\tfrac{2}{3}}{1}$ and $s_{5k+4}\in \intervalleoo{1}{+\infty}$.
\label{coro20061848}
\end{coro}

\emph{Proof}. --- Let $k\in \N^*$. As $s_{k}>0$, \eqref{eq200617234} and \eqref{eq200617235} imply that $s_{5k}\in \intervalleoo{0}{\tfrac{1}{3}}$ and $s_{5k+1}\in \intervalleoo{\tfrac{1}{3}}{\tfrac{1}{2}}$. Now \eqref{eq21060933} yields $s_{5k+2}=\frac{1+2s_{k}}{2+3s_{k}}$ whence $s_{5k+2}\in \intervalleoo{\tfrac{1}{2}}{\tfrac{2}{3}}$. However, from Proposition \ref{prop200954}, $s_{5k+3}=\frac{1}{3(1-s_{5k+2})}$. As $\tfrac{1}{2}<s_{5k+2}<\tfrac{2}{3}$, $1<3(1-s_{5k+2})<\tfrac{3}{2}$ this yields $s_{5k+3}\in \intervalleoo{\tfrac{2}{3}}{1}$. Finally, $s_{5k+4}>1$ since $s_{5k+4}=1+s_{k}$ by Corollary \ref{coro20061836}.\hfill $\square $

\bigskip

\begin{rem} --- As $s_{0}=0$, $s_{1}=\tfrac{1}{3}$, $s_{2}=1$, $s_{3}=\tfrac{2}{3}$ and $%
s_{4}=1$, we see that for every $k\in \N$, $s_{5k}\in \intervallefo{0}{\tfrac{1}{3}}$, $s_{5k+1}\in \intervallefo{\tfrac{1}{3}}{\tfrac{1}{2}}$, $s_{5k+2}\in \intervallefo{\tfrac{1}{2}}{\tfrac{2}{3}}$, $s_{5k+3}\in \intervallefo{\tfrac{2}{3}}{1}$ and $s_{5k+4}\in \intervallefo{1}{+\infty}$.
\end{rem}

\bigskip

Now we prove that the sequence $(s_{n})$ satisfies relations similar to \eqref{eq22062207} and \eqref{eq22062325}.

\begin{prop} --- For every $n\in \N^*$, $\left\lfloor s_{n-1}\right\rfloor =v_{5}(n)$. \label{prop2100949}
\end{prop}

\emph{Proof}. --- For $n=1$, $\left\lfloor s_{0}\right\rfloor =\left\lfloor 0\right\rfloor =0=v_{5}(1)$. Now assume that, for a given integer $n\geqslant 2$ and every integer $j\leqslant n-1$, $\left\lfloor s_{j-1}\right\rfloor =v_{5}(j)$. Denote $N_{k}$ ($k\in \N^*$) the parent of $N_{n}$ .

If $N_{n}$ is not the fourth child of $N_{k}$ then $5$ does not divide $n.$ Therefore $v_{5}(n)=0$ and $n-1\not\equiv 4\pmod{5}$ and, by Corollary \ref{coro20061848}, $\left\lfloor s_{n-1}\right\rfloor =0$.

If $N$ is the fourth child of $N_{k}$ then $n=5k$, whence $n-1=5k-1$. Now Corollary \ref{coro20061836} yields $s_{n-1}=1+s_{k-1},$ which implies $\left\lfloor s_{n-1}\right\rfloor =1+\left\lfloor s_{k-1}\right\rfloor $. However, by the induction hypothesis, $\left\lfloor s_{k-1}\right\rfloor =v_{5}(k)$, whence $\left\lfloor s_{n-1}\right\rfloor =1+v_{5}(k)=v_{5}(5k)$, i.e. $\left\lfloor s_{n-1}\right\rfloor =v_{5}(n)$.

Proposition 6 is proved by induction.\hfill $\square $

\bigskip

The following statement is a direct consequence of properties \ref{prop200954} and \ref{prop2100949}.

\begin{coro} --- Let $f$ be defined in \eqref{eq22062208}. Then, the sequence $(s_{n})_{n\in \N}$ satisfies $s_{0}=0$ and, for every $n\in \N^*$, 
\[s_{n}=\dfrac{1}{3(1+2\left\lfloor s_{n-1}\right\rfloor -s_{n-1})}=\dfrac{f(s_{n-1})}{3}.\]
\end{coro}

\subsection{The sequence $(s_{n})$ enumerates $\Q_{+}$}

\begin{thm} --- The mapping $n\mapsto t_{n}$ is a bijection from $\N$ to 
$\Q_{+}$. \label{th21061003}
\end{thm}

\emph{Proof}. --- As in the proof of Theorem \ref{th16061912}, we have to prove that, for every pair of non zero coprime natural integers $\left( {\alpha }\mathpunct{};{\beta }\right) $, there exists one and only one $n\geqslant 1$ such that $a_{n}=\alpha $ and $b_{n}=\beta $.

The proof is again by induction on $m=\alpha +\beta $.

If $m=2,$ then $\alpha =\beta =1$ and Corollary \ref{coro20061848} shows that $n=4$ is the only integer such that $a_{n}=b_{n}=1$.

Assume that, for a given integer $m\geqslant 2$, the property is true for every $k\in \{2,...,m\}$. Let $\left( {\alpha }\mathpunct{};{\beta }\right) $ be a pair of coprime positive integers such that $\alpha +\beta =m+1$.

As in the proof of Theorem \ref{th16061912}, we deduce from Corollary \ref{coro20061848} and Remark 2 that $n=1$ (resp. $n=2$, $n=3$ and $n=4$) if $\beta =3\alpha $ (resp. $\beta =2\alpha $, $2\beta =3\alpha $ and $\beta =\alpha $).

Now we distinguish five cases.

\emph{1st case}~: $\beta >3\alpha .$ Then, by Corollary \ref{coro20061848}, if $n$ exists, $n=5k$ with $k\in \N^*$. Hence, by \eqref{eq200617234}, $s_{k}=\frac{\alpha }{\beta -3\alpha }$. However, $\alpha +(\beta -3\alpha )=\beta -2\alpha \leqslant m$ and $\alpha $ and $\beta -3\alpha $ are coprime, which yields the conclusion by using the induction hypothesis.

\emph{2nd case}~: $2\beta <6\alpha <3\beta .$ Then, by Corollary \ref{coro20061848}, if $n$ exists, $n=5k+1$ with $k\in \N^*$. Hence, by \eqref{eq200617234}, 
\[s_{k}=\frac{3\alpha -\beta }{3\beta -6\alpha } \text{ if } 3 \nmid \beta \qquad  \text{and} \qquad s_{k} =\frac{\alpha -\frac{\beta }{3}}{\beta -2\alpha } \text{ if }3 \mid \beta\]
which yields the conclusion as in the first case. 

\emph{3rd case}~: $3\beta <6\alpha <4\beta .$ Then, by Corollary \ref{coro20061848}, if $n$ exists, $n=5k-3$ with $k\in \N^*$. Hence, by \eqref{eq200617234}, $s_{k-1}=\frac{2\alpha -\beta }{2\beta -3\alpha }$ which yields the conclusion as in the first case. 

\emph{4th case}~: $2\beta <3\alpha <3\beta .$ By Corollary \ref{coro20061848}, if $n$ exists, $n=5k-2$ with $k\in \N^*$. Then, by Proposition \ref{prop200954}, $s_{5k-1}=\frac{1}{3(1-s_{5k-2})}$ and therefore $s_{5k-2}=1-\frac{1}{3s_{5k-1}}=1-\frac{1}{3(1+s_{k-1})}$ by Corollary \ref{coro20061836}. Hence,
\[s_{k-1}=\frac{3\alpha -2\beta }{3\beta -3\alpha } \text{ if } 3 \nmid \beta \qquad  \text{and} \qquad s_{k-1} =\frac{\alpha -2\frac{\beta }{3}}{\beta -\alpha } \text{ if }3 \mid \beta\]
which yields the conclusion as in the first case. 

\emph{5th case}~: $\alpha >\beta .$ By Corollary \ref{coro20061848}, if $n$ exists, $n=5k-1$ with $k\in \N,$ $k\geqslant 2$. Then, by Corollary \ref{coro20061836}, $s_{n}=1+s_{k-1}$ and therefore $s_{k-1}=\frac{\alpha -\beta }{\beta }$ and the conclusion holds as in the first case.

Theorem 2 is proved by induction.\hfill $\square $

\section{The relation \eqref{eq02071820} with $k\geqslant 4$}

Newman result \eqref{eq22062207} and Propositions \ref{prop16061909} and \ref{prop2100949} show that the Calkin-Wilf sequence and sequences $(t_{n})$ and $(s_{n})$ are all defined by a first term $u_{0}=0$ and by a recurrence relation of the form
\[\text{for every }n\in \N^*, \quad u_{n}=\dfrac{f(u_{n-1})}{k}\]
where $k\in\{1, 2, 3\}$ and $f$ is defined by \eqref{eq22062208}. It is natural to ask if such a relation defines an enumeration of $\Q_{+}$ for every $k\geqslant 1.$ We prove now that this is
not the case.

Let $k\geqslant 4$ be an integer. Put $f_{k}=\tfrac{1}{k}f$ and consider the
sequence $(u_{n})$ defined by $u_{0}=0$ and, for every $n\in \N^*$, $u_{n}=f_{k}(u_{n-1})$. It is easy to check that the only solutions of $f_{k}(x)=x$ are
\[\gamma _{k}=\frac{1}{2}-\frac{1}{2}\sqrt{1-\frac{4}{k}}\text{ and }\delta
_{k}=\frac{1}{2}+\frac{1}{2}\sqrt{1-\frac{4}{k}},\]
and that $0<\gamma _{k}\leqslant \delta _{k}<1$. Hence, as $f_{k}$ is increasing on $\intervallefo{0}{1}$ and $f_{k}(0)=\tfrac{1}{k}>0$,  $f_{k}\left( \intervalleff{0}{\gamma_k}\right) \subset $ $\intervalleff{0}{\gamma_k}$. Moreover, $u_{1}=f_{k}(0)=\tfrac{1}{k}>u_{0}.$ Therefore $(u_{n})$ is increasing since $f_{k}$ is increasing, which proves that $(u_{n})$ is
convergent. As $f_{k}$ is continuous on $\intervalleff{0}{\gamma_k}$, $\lim u_{n}=\gamma _{k}$. Hence $\gamma _{k}$ is the only accumulation point of $(u_{n}),$ which proves that $(u_{n})$ cannot enumerate $\Q_{+},$ nor even the rationals of a given interval.

\bibliographystyle{smfalpha}
\bibliography{biblio}

\end{document}